\documentclass [12pt, a4paper,english,russian]{article}
\usepackage{amsmath, amsfonts, amssymb}
\usepackage[cp1251]{inputenc}
\usepackage[T2A]{fontenc}
\usepackage{babel,indentfirst,epigraph, cite}
\usepackage{graphicx,graphics}
\newcommand{\be}[2]{\begin{equation}\label{#1} {#2}\end{equation}}
\newcommand{\lr}[1]{\left(#1\right)}
\makeatletter
\renewcommand{\@biblabel}[1]{#1.}
\makeatother

\addto\captionsrussian{ }

\begin{document}
\break
\begin{center}
	Sitnik S.M. \\
	Inequalities for the exponential remainder of the Taylor series.\\
	Preprint, Vladivostok, 07.01.1993, 
\end{center}

\begin{center}
	\textbf{RESUME.}
\end{center}

We study sections of the exponential function Taylor series defined by \eqref{1}. Interesting inequalities for these sections were considered by G.Hardy, Kesava Menon, W. Gautschi, H.Alzer and others.   The main aim of this preprint is to investigate new proofs for the main inequality \eqref{2}  with the best possible constant and its multiple generalizations.

\newpage
\begin{center}
\large{
\textbf{COMMENTS OF 2016 year.}
}
\end{center}

\vskip 5mm

\small

From  1992 to 2016 a serious progress in this field was achieved. May be the most important progress in related to this preprint  inequalities is now for Turan -- type inequalities and serious investigation of related function classes, such as completely and absolutely monotone functions.

We add some (just some!) references related to this progress of well-known papers and some author's results.
Note that in references list at the end of the paper references of 1992 are with numbers [1]--[32] and references of 2016 are duplicated as numbers [33]--[55].

Especially we emphasize a cycle of papers of the author and Khaled Mehrez (Universite de Kairouan, Tunisia.)
Note also that in a paper [23] from the list below ([55] from the list at the end) the main inequality is not properly attributed, it was first proved by Walter Gautschi as it is duly cited in this preprint of 1992, and also in this paper author's results of this preprint and others are ignored.

\sloppy{
\begin{enumerate}

\item{\sc P. Tur\'an},
{\it On the zeros of the polynomials of Legendre},
Casopis Pest. Mat. Fys. 75 (1950), 113--122.

\item         {\sc \'A. Baricz},
        {\it Tur\'an type inequalities for hypergeometric functions},
        Proc. Amer. Math. Soc. 136 (9) (2008), 3223--3229.

\item{ \sc \'A. Baricz},
   {\it Functional inequalities involving Bessel and modified
   Bessel functions of the first kind},
    Expositiones Mathematicae, 2008, Vol. 26, P.  279--293.

\item{\sc R.W. Barnard, M.B. Gordy and K.C. Richards},
        {\it A note on Tur\'an type and mean inequalities
       for the Kummer function},
       J. Math. Anal. Appl. 349(1) (2009), 259--263.

\item{T. Pogany, Z. Tomovski. Probability distribution built by Prabhakar
function. Related Turan and Laguerre inequalities. Integral Transforms and Special Functions, 2016.}

\item {\sc  S.I. Kalmykov, D.B. Karp}, {\it Log--concavity for series in reciprocal gamma functions}, Int. Transforms and Special Func. 24 (11), 859--872, 2013.

\item{\sc  S.I. Kalmykov , D.B. Karp},
{\it  Log-convexity and log-concavity for series in gamma ratios and applications}, J. Math. Anal. Appl. 406, 400--418, 2013.

\item{\sc D.B. Karp and S.M. Sitnik},
{\it Log-convexity and log-concavity of hypergeometric-like functions},
 J. Math. Anal. Appl. 364 (2010), no. 2, 384--394.

\item{\sc  D.B. Karp ,   S.M. Sitnik,}
{\it Inequalities and monotonicity of ratios for generalized hypergeometric function},
Journal of Approximation Theory, 2009 , Vol. 161, P. 337--352.

\item        {\sc R.J. McEliece, B. Reznick and J.B. Shearer},
				{\it A Tur\'an inequality arising in information theory},
        SIAM J. Math. Anal. 12(6) (1981), 931--934.

\item        {\sc K. Mehrez, M. Ben Said, and J. El Kamel},
        {\it Tur\'an type inequalities for Dunkl and $q--$Dunkl kernel},
        arxiv.1503.04285.

\item{\sc K. Mehrez, S.M. Sitnik}, {\it Proofs of some conjectures on monotonicity of ratios of Kummer, Gauss and generalized hypergeometric functions},
arXiv:1410.6120v2 [math.CA],2014, 8 pp.

\item{\sc K. Mehrez and S. M. Sitnik},
        {\it Inequalities for sections of exponential function series and proofs of some conjectures on
        monotonicity of ratios of Kummer, Gauss and generalized hypergeometric functions},
        RGMIA Res. Rep. Collect. 17 (2014), Article ID 132.

\item {\sc K. Mehrez and S. M. Sitnik},
        {\it Proofs of some conjectures on monotonicity of ratios of Kummer and Gauss hypergeometric
        functions and related Tur\'an-type inequalities}, Analysis, 2016 (accepted for publication).

\item{\sc K. Mehrez and S.M. Sitnik},{\it On monotonicity of ratios of some hypergeometric functions}, Siberian Electronic Mathematical Reports, 2016, Vol. 13, P. 260-–268 (in Russian).

\item{\sc K. Mehrez and S.M. Sitnik},				
{\it  Monotonicity of ratios of $q$--Kummer confluent hypergeometric and $q$--hypergeometric functions and associated Tur\'an types inequalities},
 arXiv:1412.1634v1 [math.CA], 2014, 9 P.

\item{\sc K. Mehrez and S.M. Sitnik},
        {\it On monotonicity of ratios of $q$--Kummer confluent hypergeometric and $q-$hypergeometric
        functions and associated Tur\'an types inequalities}, RGMIA Res. Rep. Collect. 17 (2014), Article ID 150.

\item{\sc K. Mehrez and S.M. Sitnik},
{\it On monotonicity of ratios of some $q$--hypergeometric functions},
Matematicki Vesnik, 2016, Vol. 68, No. 3, pp. 225--231.

\item {\sc S.M. Sitnik},
{\it Inequalities for the exponential remainder},
preprint, Institute of Automation and Control Process, Far Eastern Branch of the Russian Academy of Sciences, Vladivostok, 1993 (in Russian).

\item {\sc S.M. Sitnik},
{\it A conjecture on monotonicity of a ratio of Kummer hypergeometric functions},
arXiv:1207.0936, 2012, version 2, 2014. 4 pp. (http://arxiv.org/abs/1207.0936).

\item{\sc S.M. Sitnik}, {\it  Conjectures on Monotonicity of Ratios of Kummer and Gauss Hypergeometric Functions,} RGMIA Research Report Collection, 17(2014), Article 107, 4 pp.

\item{\sc S.M. Sitnik},{\it Generalized Young and
Cauchy--Bunyakowsky Inequalities with
  Applications: A Survey}, arXiv:1012.3864, 2012, 51 pp.

\tiny{
\item{M. Merkle, Inequalities for residuals of power series: A review, Univ. Beograd, Publ. Elektrotehn. Fak. Ser. Mat 6 (1995), 79-85.}
}

\end{enumerate}
}

\newpage
\large{
\textbf{RESUME AND COMMENTS OF 1992 year.}
}

\vskip 5mm

\small

Th.1 states that \eqref{2} is equivalent to W.~Gautschi's inequality \eqref{6} from \cite{3}.

Th.2 is a simple modification, formally more general than \eqref{2}.

Th.3 contains reformations in terms of hypergeometric function, incomp\-lete gamma--function and fractional integrals.

Th.4 is a generalization of \eqref{2} again with the best constant. Proof is based on \eqref{19}--\eqref{20}, but \eqref{19} is simple and \eqref{20} is reduced to the Cebyshev inequality.

Another appropriate integral representation for the proof is $(18^*)$.

Introducing $(18^{**})$, we reduce \eqref{2} to \eqref{21}. It looks like interpolation\\ results.

Th.5 contains an interpolation inequality \eqref{22} with the best constant again. The proof is based on the Hadamard three--line theorem, \eqref{25}--\eqref{28} are consequences.

Reformulation of \eqref{2} as \eqref{29} looks like inequalities for absolutely and completely monotone functions but in different points.

Th.7 is a stronger result than \eqref{2}. The proof is reduced to an interesting inequality \eqref{32} which is trivial for integer $a=n$ but non--trivial for general $a$.

\eqref{36}--\eqref{40} are generalizations and consequences based on ideas of complete monotonicity. Kimberling's results froms \cite{22}, inequalities from \cite{19, 20, 21}.

The same results are true for more general term $R_{n, m} (x)$, occuring from the Obreshkov formula.

The interesting intriguing mistake is \eqref{42} with the wrong factorial. But if \eqref{42} were correct, then all the inequalities of this paper could be derived from the classical inequalities of Cebyshev, Jensen and Bunyakovsky by very simple procedures like \eqref{4}.

In the end there is a list of open problems and appendix with the reverse inequality \eqref{46} proved. It leads to the beautiful estimate \eqref{49} and connections with the Pade approximants.

\large{
\newpage
\noindent
	\begin{center}
		РОССИЙСКАЯ АКАДЕМИЯ НАУК
	\end{center}

	\begin{center}
		Дальневосточное отделение\\
		Институт автоматики и процессов управления
	\end{center}

\vskip 30mm
	\begin{flushright}
		Препринт
	\end{flushright}
\vskip 20mm
	\begin{center}
		Ситник С.М.
\vskip 5mm		
		Неравенства\\
		для остаточного члена ряда Тейлора\\
		экспоненциальной функции
	\end{center}
\vskip 80mm
	\begin{center}
		Владивосток\\
		1993
	\end{center}

\break
	
	\begin{flushleft}
		УДК 517.52
	\end{flushleft}

	\begin{flushleft}
		Ситник С.М. Неравенства для остаточного члена ряда Тейлора экспоненциальной функции: Препринт, Владивосток: ИАПУ ДВО РАН, 1993, 30 с.		
	\end{flushleft}
\vskip 30mm

В препринте изучаются неравенства для остатка ряда Тейлора экспоненты.

Рассматривается неравенства Хорста Алзера и Уолтера Гаучи, оценки для вырожденной гипергеометрической функции, неполной гамма-функции, дробного интеграла Римана--Лиувилля, вполне монотонных и абсолютно монотонных функций.

Работа рассчитана на специалистов по математическому Анализу.

Библ. 32.
\vskip 30mm
$\mbox{ }$\\
ОТВЕТСТВЕННЫЙ РЕДАКТОР \hfill д.ф--м.н. Быковцев Г.И. \\
РЕЦЕНЗЕНТ \hfill д.ф--м.н. Буренин А.А.~\\
Подписано к печати 19.12.1992 г.

\vskip 10mm
$\textcopyright$~~~~~~~~~~~~~~~~~~~~Институт автоматики и процессов управления

~~~~~~~~~~~~~~~~~~~~~~~ДВО РАН, 1992
}

\small
\break
	
\textbf{1.}	Рассмотрим остаточный член ряда Тейлора для экспоненциальной функции
	\be{1}{R_n (x) = \exp(x) - \sum\limits_{k=0}^n \frac{x^k}{k!} = \sum\limits_{k=n+1}^{\infty} \frac{x^k}{k!},~ x \geq 0.}

В \cite{1} Kesava Menon доказал, что справедливо неравенство
$$
R_{n-1} (x)  R_{n+1} (x) > \frac{1}{2} \lr{R_{n} (x)}^2, ~ x>0,~ n \geq 1.
$$
Недавно Horst Alzer \cite{2} доказал это неравенство с наилучшей возможной постоянной	
\be{2}{R_{n-1} (x)  R_{n+1} (x) > C_n \lr{R_{n} (x)}^2, ~ C_n = \frac{n+1}{n+2}.}
То, что постоянная $C_n$ в \eqref{2} является наилучшей, следует из
$$
\lim_{x \to 0} \frac{R_{n-1} (x)  R_{n+1} (x)}{R_{n}^2 (x)} = C_n = \frac{n+1}{n+2}.
$$
Ссылки на некоторые работы, в которых изучались свойства остатка $R_{n} (x)$, приведены в \cite{2}.

Отметим, что неравенство \eqref{2} --- это интересный и нетривиальный результат. Так лобовое применение к интегральной форме остатка
\be{3}{R_{n} (x) = \frac{1}{n!} \int\limits_0^x (x-t)^n e^t \, dt}
неравенства Буняковского по схеме
$$
\lr{\int\limits_0^x (x-t)^n e^t \, dt}^2 \leq  \lr{\int\limits_0^x \left\{ (x-t)^{\frac{n-1}{2}} e^{\frac{t}{2}} \right\} \left\{ (x-t)^{\frac{n+1}{2}} e^{\frac{t}{2}} \right\} \, dt}^2 \leq
$$
\be{4}{\leq \int\limits_0^x (x-t)^{n-1} e^t \, dt \int\limits_0^x (x-t)^{n+1} e^t \, dt}
приводит к неравенству \eqref{2} с худшей константой $n / (n+1)$, но всё--таки позволяет получить результат Kesava Menon.

Эта статья  содержит следующие результаты. Мы показываем, что неравенство \eqref{2} есть эквивалентная запись одного неравенства W.~Gau\-tschi, установленного в \cite{3}. Приведены некоторые другие эквивалентные формы неравенства \eqref{2} в терминах вырожденной гипергеометрической функции Гаусса, неполной гамма--функции,  дробного интеграла Римана--Лиувилля. Доказаны обощения неравенства \eqref{2}, указывающие на его связь с неравенством Чебышёва и теоремой Адамара о трёх прямых, неравенствами А.М.~Финка для вполне монотонных и абсолютно монотонных функций, формулой Обрешкова.

В статье используются стандартные обозначения:
$N$ --- множество натуральных чисел, $N_0$ --- неотрицательных целых чисел.

\textbf{2.} Введём, следуя \cite{3}  или \cite{4}, обозначение
$$
R_n (x) = \frac{x^{n+1}}{(n+1)!}\, e^{x Q_n (x)},~n \geq 1, ~ x \geq 0.
$$
Из теоремы Лагранжа о среднем $0< Q_n (x) <1$. Будем использовать стандартные обозначения для конечных разностей $\Delta^k$, $\nabla^k$, следуя, например, \cite{4}. Тогда основной результат из \cite{3} есть тройка неравенств
\be{5}{(-1)^k \Delta^k Q_n (x) = \nabla^k Q_n (x)>0,~k=0,1,2;~x>0.}
Рассмотрим, к каким неравенствам это приводит для величин $R_n (x)$ при $x>0$. Два первых неравенства дают эквивалентности
$$
Q_n (x) > 0 \Leftrightarrow R_n (x) > \frac{x^{n+1}}{(n+1)!}~(k=0),
$$
$$
Q_n (x) - Q_{n+1} (x) > 0 \Leftrightarrow R_n (x) > \frac{n+2}{x}\, R_{n+1} (x) ~(k=1),
$$
справедливость неравенств для $R_n$ очевидна. Рассмотрим случай $k=2$
\be{6}{Q_n (x) - 2\, Q_{n+1} (x) + Q_{n+2} (x) > 0.}

\underline{Теорема 1.} Неравенство В.~Гаучи \eqref{6} и неравенство Х.~Алзера \eqref{2} есть запись одного и того же неравенства в различных обозначениях.

Действительно, \eqref{6} эквивалентно
$$
e^{x Q_n (x)}  e^{x Q_{n+2} (x)} > \lr{e^{x Q_{n+1} (x)}}^2 \Leftrightarrow
$$
$$
\Leftrightarrow \frac{(n+1)! R_n (x)}{x^{n+1}} \frac{(n+3)! R_{n+2} (x)}{x^{n+3}} > \lr{\frac{(n+2)! R_{n+1} (x)}{x^{n+2}}}^2 \Leftrightarrow
$$
$$
R_{n} (x)  R_{n+2} (x) > \frac{n+2}{n+3}\, R_{n+1}^2 (x),
$$
но это и есть \eqref{2} с заменой $n$ на $n+1$.

Вопрос о справедливости неравенства  \eqref{6} для значений $k \geq 3$ остаётся, насколько известно автору, открытым. Он кратко обсуждается в конце статьи.

Таким образом, мы получаем следующую интерпретацию неравенства Гаучи--Алзера. Рассмотрим промежуточные значения в формуле Лагранжа для остатка ряда Тейлора экспоненты
$$
\xi_n (x) = x Q_n (x),~ R_n (x) = \frac{x^{n+1}}{(n+1)!} \, e^{\xi_n (x)} = \frac{x^{n+1}}{(n+1)!}\, e^{x Q_n (x)}.
$$
Тогда при любом $x>0$ последовательности промежуточных значений выпуклы:
$$
\nabla^2 \xi_n (x) = \Delta^2 \xi_n (x) > 0,~ \nabla^2 Q_n (x) = \Delta^2 Q_n (x) > 0.
$$

\textbf{3.} Переходя к несложным переформулировкам неравенства \eqref{2}, установим его очевидное обобщение.

\underline{Теорема 2.} Пусть $x>0$, $n,k, n-k \in N_0$. Тогда справедливо неравенство
\be{7}{R_{n-k} (x)  R_{n+k} (x) > C_{n, k} \lr{R_{n} (x)}^2}
с наилучшей возможной постоянной
\be{8}{ C_{n, k}= \frac{[(n+1)!]^2}{(n+k+1)! (n-k+1)!} = \frac{\Gamma^2 (n+2)}{\Gamma (n+2+k) \Gamma (n+2-k)}.}

Доказательство. То, что постоянная \eqref{8} наилучшая, следует из предельного равенства
$$
\lim_{x \to 0} \frac{R_{n-k} (x)  R_{n+k} (x)}{R_{n}^2 (x)} = C_{n, k}.
$$
Для доказательства  неравенства применим индукцию. При $k=1$ получаем \eqref{2}. Далее, пусть справедливы \eqref{7}--\eqref{8}. Перемножим набор неравенств \eqref{2} с различными параметрами
$$
R_{n-k-1} R_{n-k+1} > \frac{n-k+1}{n-k+2}\, R^2_{n-k},
$$
$$
R_{n-k} R_{n-k+2} > \frac{n-k+2}{n-k+3}\, R^2_{n-k+1},
$$
$$
\dots
$$
$$
R_{n+k-2} R_{n+k} > \frac{n+k}{n+k+1}\, R^2_{n+k-1},
$$
$$
R_{n+k-1} R_{n+k+1} > \frac{n+k+1}{n+k+2}\, R^2_{n+k}.
$$

В результате получим:
$$
R_{n-k-1} R_{n+k+1} > \frac{n-k+1}{n+k+2}\, R_{n-k} R_{n+k}.
$$
По предположению индукции это неравенство можно продолжить
$$
R_{n-k-1} R_{n+k+1} > \frac{(n-k+1)}{(n+k+2)}\, C_{n, k}\, R^2_{n}.
$$

Но
$$
\frac{(n-k+1)}{(n+k+2)}\, C_{n, k} = \frac{(n-k+1)}{(n+k+2)} \, \frac{[(n+1)!]^2}{(n+k+1)!\, (n-k+1)!} = C_{n, k+1}.
$$
Таким образом, получаем справедливость \eqref{7}--\eqref{8} с заменой $k$ на $k+1$, что и требовалось доказать.

Теперь сформулируем ряд утверждений, эквивалентных неравенству \eqref{7} с точной постоянной \eqref{8}.

\underline{Теорема 3.} Следующие утверждения эквивалентны \eqref{7}--\eqref{8}.

а) Пусть $n, k, n-k \in N_0$. Тогда при $x>0$ справедливо неравенство
\be{10}{ {_1F_1} (1;  n+2-k; x) {_1F_1} (1;  n+2+k; x) > \left[{_1F_1} (1;  n+2; x) \right]^2,}
причем постоянная единица в правой части является наилучшей.

б) Пусть $n, k, n-k \in N_0$. Тогда при $x>0$ справедливо неравенство
\be{11}{  \left[\gamma (n+1, x) \right]^2 < d_{n, k}\, \gamma (n+k+1, x)\, \gamma (n+1-k, x)}
с наилучшей возможной постоянной
\be{12}{ d_{n, k} = \frac{(n+1+k) (n+1-k)}{(n+1)^2} = 1 - \lr{\frac{k}{n+1}}^2.}

в) Пусть $n, k, n-k \in N_0$. Тогда при $x>0$ справедливо неравенство
\be{13}{  \left[ I^{n+1}_{0+}(e^x) \right]^2 < \frac{1}{C_{n, k}}\, I^{n+1+k}_{0+}(e^x)\, I^{n+1-k}_{0+}(e^x)}
с наилучшей возможной постоянной в правой части, причём $C_{n, k}$ определяется по формуле \eqref{8}.

В формулировке теоремы использованы обозначения: ${_1F_1} (a, b, x)$ для вырожденной гипергеометрической функции Гаусса \cite{5}, $\gamma (v, x)$ для не\-полной гамма--функции \cite{6}
\be{14} { \gamma (v, x) = \int\limits_0^x t^{v-1} e^{-t} \, dt, }
$I^{a}_{0+}  (f)$ для правостороннего интеграла Римана--Лиувилля \cite{7}.

Доказательство.

а) Эта элегантная симметричная форма неравенства \eqref{7} следует из представления
$$
R_n (x) = \frac{x^{n+1}}{(n+1)!}\,{_1F_1} (1;  n+2; x),
$$
которое может быть легко получено или найдено в \cite{8}.

б) Из \cite{6} или \cite{8} получаем связь
$$
{_1F_1} (1;  n+1; x) = \frac{n}{x^n}\, e^x \,\gamma (n, x),
$$
подстановка которой в \eqref{10} даёт
$$
\frac{(n+1-k)}{x^{n+1-k}} \, e^x \, \gamma (n+1-k, x) \, \frac{(n+1+k)}{x^{n+1+k}} \, e^x \, \gamma (n+1+k, x) >
$$
$$
> \left[ \frac{n+1}{x^{n+1}} \, e^x \, \gamma (n, x) \right]^2,
$$
что после сокращения экспонент и степеней даёт
$$
\gamma^2 (n, x) < \frac{(n+1+k) (n+1-k)}{(n+1)^2}\, \gamma (n+1+k, x)\, \gamma (n+1-k, x),
$$
а это и есть неравенство \eqref{12}.

в) По определению оператора дробного интегродифференцирования \cite{7} от функции $f(x)$
$$
I^a_{0+} (f) = \frac{1}{\Gamma (a)} \int\limits_0^x (x-t)^{a-1} f(t) \, dt.
$$
Поэтому из интегрального представления остатка \eqref{3} получаем:
\be{15} {R_n (x) = I^{n+1}_{0+} (e^x).}
Следовательно, \eqref{13} эквивалентно \eqref{7}. Теорема доказана.

Вновь отметим, что прямое применение неравенства Буняковского к определению неполной гамма--функции \eqref{14} приведёт к неравенству \eqref{11} с постоянной, равной единице, если действовать как в формуле \eqref{4}. Из \eqref{12} следует, что чем больше $k$, тем хуже постоянная единица наилучшей постоянной $d_{n, k}$, определяемой формулой \eqref{12}.

В терминах последовательностей мы получаем следующие утверждения. Если ввести обозначения
$$
a_n = \gamma (n+1, x), ~ b_n = {_1F_1} (1;  n+2; x),
$$
$$
f_n = \frac{1}{(n+1)!}  \int\limits_0^x (x-t)^n f(t) \, dt,
$$
$$
g_n (x) = (n+1)!\, R_n (x),~R_n= R_n (x),
$$
то последовательности $a_n$, $b_n$, $f_n$, $g_n$ являются логарифмически выпуклыми, а последовательность $R_n$ не является логарифмически выпуклой для всех $x>0$.  Неравенство Гаучи--Алзера эквивалентно логарифмической выпуклости последовательности $g_n$.

Отсутсвие логарифмической выпуклости для $R_n$ следует из того, что в неравенстве \eqref{2} точная  постоянная $C_n <1 $,  $(C_n)^{-1} > 1$. Вместе с тем мы можем утверждать отсутствие логарифмической  выпуклости лишь сразу при всех $x>0$,  точнее при достаточно малых $x>0$. Логически при некоторых $x$ последовательность $R_n$ может оказаться логарифмически выпуклой (сравнить с \cite{2}), хотя возможность на самом деле и не реализуется --- далее показано, что оправдано неравенство, противоположное свойству логарифмической выпуклости $R_n$.

Отметим, что для постоянной \eqref{8} также справедливо неравенство $C_{n, k}<1$ при $1 \leq k \leq n$, что эквивалентно
$$
\Gamma^2 (n+2) < \Gamma (n+2+k)\, \Gamma (n+2-k).
$$
Но последнее неравенство следует из логарифмической выпуклости \\ гамма--функции \cite{9}
$$
\Gamma^2 (x) < \Gamma (x+h)\, \Gamma (x-h),
$$
являющегося при положительных значениях всех аргументов следствием неравенства Буняковского (или, как сейчас модно говорить, неравенства КГБ).

\textbf{4.} Рассмотрим обобщения неравенства Гаучи--Алзера, основанные на неравенстве Чебышёва. При этом мы определим величину $R_n (x)$ для любых значений $a> -1$ по формуле \eqref{15}, то есть в виде дробного интеграла
\be{16}{R_a (x) = I^{a+1}_{0+} (e^x) = \frac{1}{\Gamma(a+1)}  \int\limits_0^x (x-t)^a e^t \, dt.}
Энциклопедическое изложение теории операторов дробного интегродифференцирования приведено в \cite{7}. В дальнейшем мы рассматриваем неравенства, понимая под $R_a (x)$ величину \eqref{16}.

\underline{Теорема 4.} Пусть $p>-1$, $a \geq 0$, $\beta \geq 0$.

Тогда справеделиво неравенство
\be{17}{ R_p\, R_{p+a+\beta} \geq C(p, a, \beta)\,  R_{p+a}\, R_{p+\beta} }
с наилучшей возможной постоянной
$$
C(p, a, \beta)  = \frac{\Gamma (p  + a+2) \, \Gamma (p  + \beta+2)}{\Gamma (p  +2) \, \Gamma (p+ a  + \beta+2)}.
$$

Доказательство.

Преобразуем величину $R_{a}$, интегрируя по частям,
\begin{eqnarray}
& & R_{\nu} (x) = \frac{1}{\Gamma (\nu +1)} \int\limits_0^x \left[(-1) \frac{(x-t)^{\nu+1}}{\nu+1} \right]' e^t \, dt =   \nonumber \\
& & = \frac{1}{\Gamma (\nu +2)} \left[x^{\nu+1}+ \int\limits_0^x (x-t)^{\nu+1}  e^t \, dt\right]. \label{18}
\end{eqnarray}
Представление \eqref{18} будет основой наших дальнейших выкладок. Подставляя \eqref{18} в \eqref{17}, получаем, что необходимо доказать
$$
 \frac{1}{\Gamma (p +2)} \left[x^{p+1}+ \int\limits_0^x (x-t)^{p+1}  e^t \, dt\right] \cdot
$$
$$
\cdot \frac{1}{\Gamma (p+a+\beta +2)} \left[x^{p+a+\beta+1}+ \int\limits_0^x (x-t)^{p+a+\beta+1}  e^t \, dt\right] \geq
$$
$$
\geq C(p, a, \beta)\, \frac{1}{\Gamma (p+a +2)} \left[x^{p+a+1}+ \int\limits_0^x (x-t)^{p+a+1}  e^t \, dt\right]
\cdot
$$
$$
 \cdot \frac{1}{\Gamma (p+\beta +2)} \left[x^{p+\beta+1}+ \int\limits_0^x (x-t)^{p+\beta+1}  e^t \, dt\right].
$$
Постоянные сокращаются, и мы должны доказать два неравенства
\begin{eqnarray}
& & x^{p+1} (x-t)^{p+a+\beta+1} +   x^{p+a+\beta+1} (x-t)^{p+1}  \geq  \nonumber \\
& & \geq  x^{p+a+1} (x-t)^{p+\beta+1} +   x^{p+\beta+1} (x-t)^{p+a+1},~ 0<t<x. \label{19}
\end{eqnarray}
\begin{eqnarray}
& & \int\limits_0^x (x-t)^{p+1}  e^t \, d t  \int\limits_0^x (x-t)^{p+a+\beta+1}  e^t \, dt \geq  \nonumber \\
& & \geq \int\limits_0^x (x-t)^{p+a+1}  e^t \, d t  \int\limits_0^x (x-t)^{p+\beta+1}  e^t \, dt. \label{20}
\end{eqnarray}
Для доказательства первого неравенства \eqref{19} сократим на $x^{p+1} (x-t)^{p+1}$. Получим:
$$
(x-t)^{a+\beta}  +  x^{a+\beta} \geq x^{a} (x-t)^{\beta} + x^{\beta} (x-t)^{a} , \Leftrightarrow
$$
$$
(x-t)^{\beta} [(x-t)^{a} -  x^{a}]  +  x^{\beta}[ x^{a} - (x-t)^{a}]=
$$
$$
= [x^{\beta}- (x-t)^{\beta}] [x^{a} - (x-t)^{a}] \geq 0
$$
при $0<t<x$, $a \geq 0$, $\beta \geq 0$. Второе неравенство \eqref{20} есть известное неравенство Чебышёва
$$
\int\limits_a^b p(x)\, dx  \int\limits_a^b p(x) f(x) g(x )\, dx  \geq \int\limits_a^b p(x) f(x)\, dx \int\limits_a^b p(x) g(x)\, dx
$$
для $p(x) \geq 0$, $f$ и $g$ одинаково монотонных (см. \cite{10}--\cite{14}). В нашем случае $a=0$, $b=x$, $p(x) = (x-t)^{p+1} e^t \geq 0$, $f(x) = (x-t)^{a}$, $g(x) = (x-t)^{\beta}$, причём $f$ и $g$ одновременно убывают при $a \geq 0$, $\beta \geq 0$. Следовательно, справедливо \eqref{20}, а тогда и \eqref{17}. Отметим, что следуя \cite{19}, можно свести \eqref{20} и к интегральному неравенству о средних.

Неулучшаемость постоянной видна при $x \to 0$, что завершает доказательство теоремы 4.

Частными случаями неравенства \eqref{17} при $a=\beta=1$ являются неравенство Гаучи--Алзера \eqref{2} (в котором $n$ заменено на $n+1$); при $p=n-k$, $a=\beta=k$ неравенство \eqref{7}, а также неравенство \eqref{9}.

Возможно реализовать ту же идею на основе другого интегрального представления. Для $n \in N$ оно получается так:
$$
R_n (x) = \frac{x^{n+1}}{(n+1)!}\, {_1F_1} (1; n+2; x) = \frac{x^{n+1}}{(n+1)!}\, e^x \, {_1F_1} (n+1; n+2; -x) =
$$
$$
= \frac{x^{n+1}}{n!} \, e^x  \int\limits_0^1 t^n e^{-xt} \, dt =  \frac{x^{n+1}\, e^x}{(n+1)!}  \lr{e^{-x} + x \int\limits_0^1 t^{n+1} e^{-xt} \, dt}=
$$
$$
=\frac{x^{n+1} }{(n+1)!}  \lr{1 + x  e^x\int\limits_0^1 t^{n+1} e^{-xt} \, dt}. \eqno{(18^*)} \label{18*}
$$

Теперь, например, неравенство Гаучи--Алзера сводится к проверке двух неравенств:
$$
t^n + t^{n+2} \geq 2 \, t^{n+1} \Leftrightarrow t + \frac{1}{t} \geq 2,~t>0,
$$
$$
\int\limits_0^1 t^{n} e^{-xt} \, dt \int\limits_0^1 t^{n+2} e^{-xt} \, dt \geq \lr{\int\limits_0^1 t^{n+1} e^{-xt} \, dt}^2.
$$
Но первое неравенство очевидно, а второе следует из неравенства Буняковского или интегрального неравенства о средних.

\textbf{5.} Запишем неравество Гаучи--Алзера \eqref{2} в виде
\be{21}{B_n^2 \leq B_{n-1} B_{n+1},}
где по определению полагаем:
$$
B_{\nu} = \Gamma (\nu+2)\, R_{\nu} (x).
\label{18**}
\eqno{(18^{**})}
$$
Форма \eqref{21} напоминает интерполяционные неравенства \cite{15}. Поэтому правдоподобно, что для доказательств можно использовать теорему Адамара о трёх прямых \cite{15}.

\underline{Теорема 5.} Пусть $\nu > -1$, $a \geq 0$, $0 \leq \theta \leq  1$, $x>0$. Тогда справедливо интерполяционное неравенство с наилучшей постоянной
\be{22}
{
	\begin{array}{l}
		R_{\nu + a \theta} (x) \leq C (\nu, a, \theta) \left\{R_{\nu} (x) \right\}^{1-\theta} \left\{R_{\nu+a} (x) \right\}^{\theta}, \\
		C (\nu, a, \theta)  = \frac{ \left\{\Gamma (\nu+2) \right\}^{1-\theta} \left\{\Gamma (\nu+a+2) \right\}^{\theta}}{\Gamma (\nu+a \theta+2)}
	\end{array}	
}

Доказательство. В силу \eqref{18}
$$
B_{\nu} (x)\, x^{\nu +1 } + \int\limits_0^x (x-t)^{\nu +1} e^t \, dt.
$$
Введём в расмотрение функцию от комплексной переменной $S$
$$
f (S) = x^{\nu +1 + a S} +  \int\limits_0^x (x-t)^{\nu +1 + a S} e^t \, dt, ~ a>0.
$$
Функция $f(S)$ является аналитической в полосе $0 < Re \, s < 1$ при $x>0$. Кроме того, она ограничена на всех прямых $Re \, S = \theta$, $0 \leq \theta \leq 1$,
\begin{eqnarray}
& &  \sup\limits_{\substack{S=i \xi + \theta \\ \xi \in R}} |f(S)| \leq \sup\limits_{\substack{S=i \xi + \theta \\ \xi \in R}} \lr{ |x^{\nu +1 + a \theta + i a \xi}| + \int\limits_0^x |(x-t)^{\nu +1 + a \theta + i a \xi}| e^t \, dt} = \nonumber \\
& & =  x^{\nu +1 + a \theta } + \int\limits_0^x (x-t)^{\nu +1 + a \theta } e^t \, dt. \label{23}
\end{eqnarray}
Следовательно,
$$
\sup\limits_{\substack{S=i \xi + \theta \\ \xi \in R}} |f(S)| = B_{\nu+a \theta},
$$
так как супремум в \eqref{23} достигается при $\xi = 0$. Итак, $f(S)$ ограничена и непрерывна при $Re \, S =0$, $Re \, S = 1$, и по теореме Адамара о трёх прямых получаем неравенство
\be{24}{B_{\nu+a \theta} \leq B_{\nu}^{1-\theta} B_{\nu+a}^{\theta},}
что эквивалентно
$$
\Gamma (\nu + a \theta + 2)  R_{\nu + a \theta} (x) \leq
$$
$$
 \leq \left\{ \Gamma (\nu + 2) \, R_{\nu } (x)\right\}^{1- \theta} \left\{ \Gamma (\nu + a + 2) \, R_{\nu + a} (x)\right\}^{ \theta}.
$$
Для завершения доказательства осталось проверить точность постоянной, перейдя к пределу $x \to 0$.

Анализ доказательства показывавет, что важным моментом является переход к величинам $B_{\nu}$, интегральное представление которых не содержит зависящих от $\nu$ постоянных.

Частным случаем \eqref{22} является неравенство \eqref{7}, которое получается при выборе параметров $\theta =1/2$, $\nu = n - k$, $a = 2 k$.

\underline{Следствие.} Справедливы неравенства:

а) при $p \geq 1$, $\nu > -1$, $a > 0$
\be{25}{ R^p_{\nu + \frac{a}{p}} (x) \leq C_1 (\nu, a, p) R_{\nu+a} (x) R^{p-1}_{\nu} (x),}
$$
C_1 (\nu, a, p) = \frac{[\Gamma(\nu+2)]^{p-1} \Gamma(\nu+a+2)}{[\Gamma(\nu+\frac{a}{p}+2)]^{p}};
$$

б) при $n, k \in N_0$
\be{26}{R_{n+1}^k (x) \leq \frac{(n+k+1)!}{(n+2)^{k-1}\, (n+2)!} \,  R_{n+k} (x) \, R_{n}^{k-1} (x);}

в) при $n, k \in N_0$, $a \geq \beta \geq 0$
\be{27}{R_{n+\beta}^{a} (x) \leq C_2 (n, a, \beta) \,  R_{n}^{a - \beta} (x) \, R_{n+a}^{\beta} (x),}
$$
C_2 (n, a, \beta)  = \lr{\frac{\Gamma(n+2)}{\Gamma(n+\beta+2)}}^{a} \lr{\frac{\Gamma(n+a+2)}{\Gamma(n+2)}}^{\beta}.
$$
Каждое из неравенств \eqref{25}--\eqref{27} влечёт \eqref{2}. Например, в \eqref{26} для этого нужно выбрать $k=2$ и заменить $n+1$ на $n$. Все постоянные в неравенствах \eqref{25}--\eqref{27} неулучшаемы.

\underline{Теорема 6.} Пусть $\nu>-1$, $a \geq 0$, $0 \leq \theta_i \leq 1$, $1 \leq i \leq k$, $\theta_1+\theta_2+\dots+\theta_k = \beta$, $x>0$. Тогда
\be{28}{\prod_{i=1}^{k} R_{\nu + a \theta_i} (x) \leq C(\theta) \, R_{\nu}^{k-\beta} (x) \, R_{\nu+a}^{\beta} (x),}
где $C(\theta) = \prod\limits_{i=1}^{k} C(\nu, a, \theta_i)$ --- есть произведение постоянных из \eqref{22}.

\textbf{6.} В этом пункте мы рассмотрим обощения неравенства Гаучи--Ал\-зера \eqref{2} с использованием идей, связанных с абсолютно и вполне монотонными функциями. Наводящим соображением является перезапись \eqref{2} для $f(x) = \exp(x)$ в виде
\be{29}{f^{(n+1)} (\xi_1 (x)) \cdot  f^{(n-1)} (\xi_2 (x)) \geq \left[ f^{(n)} (\xi_3 (x))\right]^2. }
Это напоминает неравенство между производными абсолютно или впол\-не монотонных функций, только производные вычисляются в различных точках. Теория абсолютно или вполне монотонных функций, восходящая ещё к Хаусдорфу и Бернштейну, изложена у Уиддера в \cite{16} (см. также \cite{17, 18, 19, 20, 21, 22, 23}).

Неравенства типа \eqref{29} (при одинаковых значениях аргументов) установлены в \cite{16}. Они обощены А.М.~Финком в \cite{21} на произведения любого числа производных. Дальнейшие обобщения и общий элегантный метод получения подобных неравенств, сводящий их к неравенствам Чебышёва, предложен в \cite{19}--\cite{20}.

Преобразуем величину $R_{\nu}(x)$
$$
 R_{\nu}(x) = \frac{1}{\Gamma(\nu+1)} \int\limits_0^x (x-t)^{\nu} e^t \, dt = \frac{1}{\Gamma(\nu+1)} \int\limits_0^x y^{\nu} e^{x-y} \, dy =
$$
$$
 = \frac{e^x\, x^{\nu+1}}{\Gamma(\nu+1)} \int\limits_0^x t^{\nu} e^{-xt} \, dt = \frac{e^x\, x^{\nu+1}}{\Gamma(\nu+2)} \int\limits_0^{\infty} e^{-xt} \, d w(t) , \label{18***} \eqno{(18^{***})}
$$
где неубывающая функция $w(t)$ равна
$$
w(t) = \left\{
\begin{array}{ll}
 t^{\nu+1}, & 0 \leq t \leq 1 \\
 1, & t>1.
\end{array}
  \right.
$$

\underline{Определение.} (\cite{16, 19}). Функция, представимая в виде преобразования Лапласа
$$
f(x) = \int\limits_0^{\infty} e^{-xt} \, d w(t)
$$
с непрерывной неубывающей ограниченной функцией $w(t)$, называется вполне монотонной.

Доказано, что это определение эквивалентно выполнению трёх условий \cite{23}, \cite{16}:

1. $f(x) \in C^{\infty} (0, \infty)$.

2. $(-1)^k f^{(k)} (x) \geq 0$, $x \geq 0$.

3. $f(x) = o (1)$ при $x \to \infty$.\\
Для вполне монотонных функций справедливы неравенства между производными, простейшее из которых имеет вид \cite{19}
\be{30}{f''(x) f(x) \geq [f'(x)]^2.}

Перейдём к доказательству неравенства, которое уточняет \eqref{2}, позволяя оценить разность в \eqref{2} между левой и правой частями.

\underline{Теорема 7.} Пусть $a > 0$, а величина $R_{a} (x)$ определена по формуле \eqref{16}. Тогда справедливо неравенство
\be{31}{R_{a+1} (x) R_{a-1} (x) - \frac{a+1}{a+2} \, R_{a}^2 (x) \geq \frac{1}{a+2} \, R_{a}^2 (x) - \frac{a+2}{x^2} \, R_{a+1}^2 (x) \geq 0. }

Доказательство. Рассмотрим функцию
$$
f(x) = e^{-x} x^{- (\nu+1)} R_{\nu} (x),~\nu>1.
$$
Мы уже показали, что эта функция является вполне монотонной, то есть представима в виде преобразования Лапласа с нужной функцией $w(t)$. Так как
$$
R_{\nu} (x) = I_{0+}^{\nu+1} (e^x),
$$
то справедливы формулы
$$
R'_{\nu} (x)=R_{\nu-1} (x),~R''_{\nu} (x) = R_{\nu-2} (x),
$$
которые можно получить непосредственно. Они являются частными случаями полугруппового свойства операторов дробного интегродифференцирования \cite{7}
$$
\frac{d}{d x} I^{\beta}_{0+} = I^{\beta-1}_{0+},~ \beta>0.
$$
Вычисляем производные функции $f(x)$
$$
f'(x) = e^{-x} x^{-(\nu+1)} \lr{R_{\nu-1} - \lr{1+\frac{\nu+1}{x}} R_{\nu}},
$$
$$
f''(x) = e^{-x} x^{-(\nu+1)} \left[ R_{\nu-2} - 2 \lr{1+\frac{\nu+1}{x}} R_{\nu-1} \right. +
$$
$$
+ \left. \lr{1 + \frac{2(\nu+1)}{x} +\frac{(\nu+1)(\nu+2)}{x^2} } R_{\nu} (x) \right].
$$
После подставновки в \eqref{30} и скучных вычислений получаем \eqref{31}, в котором обозначено $\nu-1=a$. Осталось проверить неравенство
\be{32}{R_{a} (x) \geq \frac{a+2}{x} R_{a+1}.}
Для натуральных $a=n$ это очевидно и уже встречалось нам как одно из неравенств \eqref{5} при $k=1$. При остальных $a$ необходимо доказать неравенство
$$
\frac{1}{\Gamma(a+1)} \int\limits_0^x (x-t)^{a} e^t \, dt \geq \frac{(a+2)}{x} \frac{1}{\Gamma(a+2)} \int\limits_0^x (x-t)^{a+1} e^t \, dt.
$$
Оно эквивалентно соотношению
$$
\int\limits_{0}^x \left[(x-t)^a - \frac{a+2}{a+1} \, \frac{(x-t)^{a+1}}{x} \right] e^t \, dt \geq 0,
$$
которое не является тривиальным, так как функция под знаком интеграла меняет знак. Выполнив замену переменных по формуле $t=x y$, получим:
$$
x^{a+1} \int\limits_0^1 \left[ (1-y)^{a} - \frac{a+2}{a+1} (1-y)^{a+1} \right] e^{xy} \, dy = x^{a+1} \int\limits_0^1 k(y)  \varphi_x (y) \, dy.
$$

Далее мы используем результат из \cite{24} (глава X).

\underline{Лемма.} Для того, чтобы интегрируемая функция $k(y)$ обладала свойством
$$
\int\limits_0^1 k(y)  \varphi (y) \, dy \geq 0
$$
для всех положительных, возрастающих и ограниченных $\varphi (y)$, необходимо и достаточно, чтобы выполнялось условие
\be{33}{\int\limits_0^1 k(y)  \, dy \geq 0,~ 0 \leq p \leq 1.}

Очевидно, что функция $\varphi_x (y) = \exp (xy)$ при каждом фиксированном $x$ удовлетворяет условию леммы, поэтому осталось проверить условие \eqref{33}. Получаем, обозначив $1-p=q$,
$$
\int\limits_p^1 \left[ (1-y)^{a} - \frac{a+2}{a+1} (1-y)^{a+1} \right] dy =
$$
$$
= \int\limits_0^q \left[ t^{a} - \frac{a+2}{a+1} \, t^{a+1} \right] d t = \frac{1}{a+1}\, q^{a+1} (1-q) \geq 0,
$$
так как по условию $0 \leq q \leq 1$. Теорема доказана.

Отметим, что встретившееся неравенство
\be{34}{I^{a}_{0+} (f) \geq \frac{a+1}{x} I^{a+1}_{0+} (f)}
для положительной, возрастающей и ограниченной функции $f(x)$ представляет, на наш взгляд, определённый самостоятельный интерес.

Использование более общих неравенств А.М.~Финка \cite{21}, а также их обобщений  с использованием идей выпуклости и мажоризации \cite{19}--\cite{20} приводит к другим неравенствам для величин $R_{\nu} (x)$. Рассмотрим лишь простейшие из них, вытекающие из того, что введённая при доказательстве теоремы 7 функция является вполне монотонной.

Так, использование свойства $f'(x) \leq 0$ даёт в совокупности с \eqref{32} двустороннее неравенство
\be{35}{\frac{\nu+1}{x}\, R_{\nu} (x) \leq R_{\nu-1} (x) \leq \lr{1+\frac{\nu+1}{x}} R_{\nu} (x).}

Поэтому, если ввести обозначение для разности
$$
\varepsilon_{\nu} (x) = \frac{R_{\nu} (x)}{R_{\nu+1} (x)} - \frac{(\nu+2)}{x},
$$
то из \eqref{35} следует, что $ 0 \leq \varepsilon_{\nu} (x) \leq 1$, причём крайние  значения достигаются при $x \to 0$ и $x \to \infty$ соответственно. Можно утверждать, что при фиксированном $\nu$ функция $\varepsilon_{\nu} (x)$ изменяется между этими значениями, монотонно возрастая. Действительно, $\varepsilon'_{\nu} (x) \geq 0$ эквивалентно \eqref{31}.

Неравенство $f''(x) \geq 0$ приводит к
$$
R_{\nu-2} (x) - 2 \lr{1+\frac{\nu+1}{x}} R_{\nu-1} + \lr{1 + \frac{2(\nu+1)}{x} +\frac{(\nu+1)(\nu+2)}{x^2} } R_{\nu} \geq 0.
$$
С другой стороны, согласно \cite{19}, из $(a-1)^2 \geq 0$ следует, что
$$
f''(x) + 2 f'(x) + f(x) \geq 0.
$$
Это приводит к усилению неравенства \eqref{32}
\be{36}{R_{\nu-2} (x) - \frac{\nu}{x}\, R_{\nu-1} (x) \geq \frac{\nu+2}{x} \lr{  R_{\nu-1} (x)- -\frac{\nu+1}{x} \, R_{\nu} (x)}  \geq 0.
	}

Интересны следствием интегрального представления ($18^{***}$) для \\ вполне монотонной функции является возможность применения неравенств Кимберлинга \cite{22}. Из ($18^{***}$) следует, что $g: [0, \infty) \to (0, 1]$, где
$$
g(x) = \Gamma (\nu +2) e^{-x} x^{- (\nu+1)} R_{\nu} (x) = (\nu+1) \int\limits_0^1 t^{\nu} e^{-xt} \, dt \leq 1.
$$

Поэтому из \cite{22}, \cite{19} следуют неравенства
\be{37}{R_{\nu} (x+y) \geq \Gamma (\nu+2) \lr{\frac{1}{x}+\frac{1}{y}}^{\nu+1} R_{\nu} (x)\, R_{\nu} (y),}
\be{38}{R_{\nu} (x+y) \leq  \left[ \frac{x+y}{(x+py)^{\frac{1}{p}} x^{\frac{1}{q}}} \right]^{\nu+1} R_{\nu}^{\frac{1}{p}} (x+py)\, R_{\nu}^{\frac{1}{q}} (x),}
где $1/p+1/q = 1$,
\be{39}{R_{\nu} (2 x) \geq \frac{\Gamma (\nu+2) \, 2^{\nu+1}}{x^{\nu+1}}\, R^2_{\nu} (x),}
\be{40}{R_{\nu} (x+y) \leq \left[ \frac{(x+y)}{4xy}\right]^{\nu+1} R_{\nu} (2x)\, R_{\nu} (2y).}
Отметим, что неравенство \eqref{40} точное при $x \to y$, а \eqref{39} --- при $x \to 0$.

\textbf{7.} Сделаем ряд замечаний.

Исследование остаточного члена
\be{41}{R_n (-x) = e^{-x} - \sum\limits_{k=0}^n \frac{(-x)^k}{k!} = \frac{(-1)^{n+1}}{n!} \int\limits_0^x (x-t)^n e^{-t}\, dt}
в силу приведённого интегрального представления позволяет использовать ту же технику и получить по существу те же неравенства \eqref{2}, \eqref{7}, \eqref{9}, \eqref{17}, \eqref{22}, \eqref{25}--\eqref{27}, \eqref{31} для модуля величины \eqref{41}. При этом доказательство теоремы 7 потребует рассмотрения абсолютно монотонных функций (вместо вполне монотонных) и неравенств А.М.~Финка для них.

Величина $R_n (x)$ входит в решение следующей интересной задачи: для уравнения
$$
y'(x) = \lambda y(x),~y(0) = y_0
$$
ошибка при решении этого уравнения любым $n$ --- шаговым методом Рунге--Кутта равна
$$
\varepsilon_n = R_n (\lambda h) \cdot y_0,
$$
где $h$ -- шаг метода. Это доказано в \cite{25} (глава 2). Таким образом, все неравенства для величин $R_n$ приводят к неравенствам для ошибок $\varepsilon_n$. Этот результат переносится на явные $n$--шаговые методы решения абстрактного уравнения $u'=A u$.

Возможное обобщение наших оценок связано с формулой Обрешкова, которая является обощением формулы Тейлора и служит основой для вывода аппроксимаций Паде экспоненты \cite{26, 27} (экспаденты по терминологии \cite{27}, \cite{32}).
$$
\sum\limits_{j=0}^m (-1)^j \frac{C^j_m}{C^j_{n+m}} \frac{(x-x_0)^j}{j!} f^{(j)} (x) =
$$
$$
= \sum\limits_{j=0}^n (-1)^j \frac{C^j_n}{C^j_{n+m}} \frac{(x-x_0)^j}{j!} f^{(j)} (x_0)+
$$
$$
+ \frac{1}{(n+m)!} \int\limits_{x_0}^x (x-t)^n (x_0-t)^m f^{(n+m+1)} (t)\, dt.
$$
Поэтому для величин
$$
R_{n, m} (x) = \lr{\sum\limits_{j=0}^m (-1)^j \frac{C^j_m}{C^j_{n+m}} \frac{x^j}{j!}} e^x - \sum\limits_{j=0}^n  \frac{C^j_n}{C^j_{n+m}} \frac{x^j}{j!} =
$$
$$
= \frac{(-1)^m}{(n+m)!} \int\limits_0^x (x-t)^n t^m e^t \, dt,~R_{n, 0} (x) = R_n (x)
$$
можно вывести аналоги всех неравенств настоящей работы. Они сводятся к нашим при $m=0$.

В частности, при $m=1$ получаем:
$$
R_{n, 1} (x) = \lr{1- \frac{1}{n+1} x} e^x - \sum\limits_{j=0}^n \lr{1 - \frac{j}{n+1}\, x}  \frac{x^j}{j!}.
$$
Это выражение содержит средние Чезаро \cite{27, 29, 30} для частичной суммы Тейлора экспоненты
$$
\gamma_n^{(1)} (x) = \sum\limits_{j=0}^n \lr{1 - \frac{j}{n+1}\, x} \frac{x^j}{j!}.
$$
Таким образом, можно получить неравенства, в которые входят  величины $\gamma_n^{(1)} (x)$ --- средние
Чезаро. О предельных соотношениях для этих средних см. \cite{27}, \cite{32}.

Интересно отметить, что все неравенства этой работы могут быть получены из одной очень плодотворной ошибки. Если взять формулу \eqref{3} с неверным факториалом
\be{42}{R_n (x) \stackrel{?}{ =} \frac{1}{(n+1)!} \int\limits_0^x (x-t)^n e^t \, dt,}
то простое применение неравенства Буняковского по схеме \eqref{4} приводит к неравенству Гаучи--Алзера \eqref{2} и его обобщению \eqref{7}, применение неравенства Чебышёва по той же схеме даёт \eqref{17}, обобщенного неравенства Чебышёва или неравенства Иенсена --- даёт \eqref{25}, неравенства для интегральных средних --- даёт \eqref{27}. Таким образом, все выводы, следующие из этой ошибки  и применения классических неравенств Анализа, оказываются верными и доказаны в настоящей работе. Более того, именно неверная формула \eqref{42}, из которой всё так просто выводится, стимулирует использование интегральных представлений \eqref{18}--$(18^{*})$ -- $(18^{**})$ -- $(18^{***})$, основная цель которых состоит в замене постоянной $n!$ на $(n+1)!$

Поэтому символ (!) в конце предыдущей строки автор рассматривает не только как символ факториала, но и как восклицательный знак --- выражение своего удивления и благодарности по отношению к ошибке \eqref{42}.

В заключение приведём несколько задач, решения которых неизвестны автору.

1. Верно ли, что функция
$$
f(x) = \frac{R_{n-1} (x) R_{n+1} (x)}{R_{n}^2 (x)}
$$
является монотонно возрастающей по $x$ на $(0, \infty)$? Отметим, что неравенство \eqref{2} имеет вид $f(x) > f(0)$, $x>0$. Дифференцирование сводит дело к проверке неравенства
$$
R_{n-2}\, R_{n} \, R_{n+1} + R_{n-1}\, R_{n}^2 > 2 R_{n-1}^2 \, R_{n+1}.
$$

2. Для каких ещё гипергеометрических функций и параметров в них справедливо \eqref{10}?

По--видимому, этот вопрос можно связать с абсолютно монотонными функциями.

3. Справедливо ли неравенство \eqref{2} с некоторой постоянной в другую сторону

4. Можно ли выводить неравенства типа \eqref{2} из разложения $R_{n} (x)$ в цепную дробь \cite{28}? В применении к неравенствам для функций Бесселя и Макдональда этот подход через цепные дроби реализованы в \cite{27}, \cite{32}.

5. Являются ли диагональные аппроксимации Паде экспоненты абсолютно монотонными при $x > 0$ и вполне монотонными при $x<0$? К каким неравенствам это приводит?

6. Является ли последовательность $R_n (x)$, $x>1$ логарифмически выпуклой хотя бы при одном $x$?

7. Верно ли, что
$$
	\lim\limits_{n \to \infty} \lr{ \frac{x}{n} \frac{R_{n-1}(x) }{R_n (x)} } = e^{m+1}?
$$

8. Справедливы ли неравенства Гаучи \eqref{5} для каких--нибудь $k \geq 3$? При $k=3$ соответствующее неравенство для $R_n (x)$ принимает вид
$$
R_n \, R_{n+2}^3 > \frac{(n+2)(n+4)}{(n+3)^2} \, R_{n+1}^3 \, R_{n+3}.
$$

9. Вычислить предел (если он существует)
$$
\lim\limits_{n \to \infty} \frac{R_{n, m} (a n)}{e^{a m}} = L (a, m).
$$
О случае $m=0$ см. \cite{28}.

\break
\begin{center}
\underline{Дополнение}
\end{center}

В этом дополнении приводятся решения задачи 3,6 со стр. 20--21. Полученный результат эквивалентен неравенствам между экспонентой и её аппроксимациями первой строки таблицы Паде.

\underline{Теорема 8.} Пусть $n \in N$. Тогда справедливо обращение неравенства Гаучи--Алзера
\be{43}{R_n^2 (x) > R_{n-1 } (x)\, R_{n+1 } (x).}
Постоянная единица в правой части является наилучшей и достигается при $x \to \infty$.

Доказательство. Напомним, что
$$
R_{n-1 } (x) = \frac{x^n}{n!} + \dots,~R_{n} (x) = \frac{x^{n+1}}{(n+1)!} + \dots,~R_{n+1 } (x) = \frac{x^{n+2}}{(n+2)!} + \dots,
$$
$$
R_{n-1 } (x) = R_{n} (x) + \frac{x^n}{n!},~R_{n+1 } (x) =R_{n} (x)-\frac{x^{n+1}}{(n+1)!}.
$$
Выражая в \eqref{43} все величины через $R_n$, получаем:
$$
R_{n} (x) \lr{\frac{x^n}{n!} - \frac{x^{n+1}}{(n+1)!} } < \frac{x^{2n+1}}{n!\, (n+1)!},
$$
или после сокращений
\be{44}{(n+1-x)\, R_n(x) < \frac{x^{n+1}}{n!}.}
Последнее неравенство эквивалентно
$$
(n+1-x) \lr{\frac{1}{n+1} + \frac{x}{(n+1)(n+2)} + \frac{x^2}{(n+1)(n+2)(n+3)} +\dots} < 1.
$$
Преобразуем левую часть
$$
1 - \frac{x}{n+1} + \frac{x}{n+2} - \frac{x^2}{(n+1)(n+2)} + \frac{x^2}{(n+2)(n+3)} - \dots =
$$
$$
 = 1 - x \lr{ \frac{1}{n+1} - \frac{1}{n+2}} -x^2 \lr{\frac{1}{(n+1)(n+2)} - \frac{1}{(n+2)(n+3)}} - \dots =
$$
$$
= 1 - \frac{x}{(n+1)(n+2)}  - \frac{2 x^3}{(n+1)(n+2)(n+3)}-
$$
$$
 - \frac{3 x^3}{(n+1)(n+2)(n+3) (n+4)} - \dots < 1
$$
при $x>0$, что очевидно. Теорема доказана.

Преобразуем \eqref{44}к виду неравенств для самой экспоненты. Тогда для первых $n$ получаем при $x>0$

1) $n=0$,
$$
e^x < \frac{1}{1-x}, ~~~~~~x<1,
$$

2) $n=1$,
$$
e^x < \frac{2+x}{2-x}, ~~~~~~x<2,
$$

3) $n=3$,
$$
e^x < \frac{x^2 +4 x +6}{6-2x}, ~~~~~~x<3,
$$
и обратные неравенства при $x > n+1$, $n = 0, 1, 2$. В правой части выделились элементы первой строки таблицы Паде экспоненты, $[0/1]$, $[1/1]$ и $[2/1]$. Такими наводящими соображениями получается

\underline{Теорема 9.} Пусть $R_{n, 1}$ --- экспаденты $[n/1]$, $n \geq 0$. Тогда справедливо неравенство
\be{45}{e^x < R_{n, 1} (x),~x<n+1,}
а при $x > n+1$ справедливо обратное неравенство.

Можно прямо вывести \eqref{45} из \eqref{44} с использованием явного вида $R_{n, 1} (x)$ \cite{31, 26, 27} . Мы приведём несколько другое доказательство, которое выявляет связь \eqref{45} и \eqref{43} и не использует явного вида аппроксимаций $R_{n, 1} (x)$.

Пусть $f(x)$ --- функция , представимая при $x \geq 0$ рядом Тейлора
$$
f(x) = t_n(x) + R_n (x),~t_n (x) = \sum_{k=0}^{n} \frac{f^{(k)} (0)}{k!}\, x^k,
$$
$$
\Delta (x) = x f^{(n+1)} (0) - (n+1) f^{(n)} (0),
$$
$R_{n, 1}(x)$ --- элементы первой строки таблицы Паде для функции $f(x)$.

\underline{Теорема 10.} Неравенство для остатков ряда Тейлора функции $f(x)$
\be{46}{R_n^2 (x) > R_{n-1 } (x)\, R_{n+1 } (x)}
эквивалентно неравенству
\be{47}{f(x) > R_{n, 1} (x)}
при $\Delta (x) > 0$ и эквивалентно обратному неравенству при $\Delta (x) < 0$.

\underline{Доказательство.} После подстановок \eqref{46} принимает вид
$$
f(x) \lr{t_{n+1} (x) + t_{n-1} (x) - 2\, t_n (x)} > t_{n-1} (x) t_{n+1} (x) - t_n^2(x).
$$
Формула Тейлора даёт
$$
t_{n+1} (x) + t_{n-1} (x) - 2\, t_n (x) = \frac{x^{n+1}}{(n+1)!} \, f^{(n+1)} (0) -
$$
$$
- \frac{x^{n}}{n!}\, f^{(n)} (0) = \frac{x^{n}}{(n+1)!}\, \Delta (x).
$$
Поэтому при $\Delta (x) > 0$ получаем:
\be{48}{f(x) > \frac{t_{n-1} (x)\, t_{n+1} (x) - t_n^2(x)}{t_{n+1} (x) + t_{n-1} (x) - 2\, t_n (x)}.}
Правая часть \eqref{48} --- это преобразование Эйткена последовательности $t_n$ \cite{28, 31}. Частичные суммы ряда Тейлора --- это нулевая строка таблицы Паде. Но известно, что преобразование Эйткена нулевой строки есть в точности первая строка таблицы Паде \cite{31}. Следовательно, \eqref{48} эквивалентно \eqref{47}. При $\Delta (x) < 0$ получаем неравенство со знаком ``$<$''. Теорема доказана.

Чтобы вывести теорему 9 из теоремы 10 достаточно заметить, что для экспоненты $\Delta (x) = x - (n+1)$.

Таким образом, свойство противоположное логарифмической выпуклости последовательности остатков ряда Тейлора \eqref{46} эквивалентно не\-равенствам между самой функцией и её аппроксимациями Паде $[n/1]$.

\underline{Следствие 1.} Справедливо двустороннее неравенство с точными константами
\be{49}{\frac{n+1}{n+2} < \frac{R_{n-1}(x)\, R_{n+1}(x)}{R^2_n (x)} <1.}

\underline{Следствие 2.} Вопрос 6 на стр. 21 имеет отрицательный ответ.

Приведём ещё одно неожиданное следствие  неравенств, доказанных в этом дополнении. Рассмотрим стандартные разностные схемы для уравнения теплопроводности: явную, неявную и Кранка--Никольсона. Тогда при выполнении соответствующих естественных условий приближения пр помощи этих схем заключают точное решение в вилку. А именно, приближение по явной схеме всегда меньше точного решения, а приближения по неявной схеме или Кранка--Никольсона всегда больше точного решения.

Дополним список задач на стр. 21--22.

10. Пусть $R_{n, m} (x)$ --- экспадента \cite{27} вида $[n/m]$, $a_1$ --- первый положительный корень знаменателя $Q_m$ экспаденты $R_{n, m}$. Тогда
$$
(-1)^m e^m > (-1)^m R_{n,m}(x),~0<x<a_1,
$$
а после каждого корня $a_i$ знак неравенства меняется на противоположный.

При $m=0$ это очевидно, при $m=1$ доказано в настоящей работе.

Замечание к задаче (1) на стр. 20.

Неравенство \eqref{49} делает утверждение о монотонности функции $f(x)$ ещё более правдоподобным. Теперь \eqref{49} позволяет получить
$$
R_{n-2}\, R_n\, R_{n+1}  + R_{n-1}\, R^2_n > \lr{1+\frac{n}{n+1}}\, R^2_{n-1}\, R_{n+1},
$$
что близко к требуемому неравенству, но слабее его.

11. Определим функцию при $x>0$.
$$
g_n (x) = \frac{R_{n-1} (x)}{R_n (x)}.
$$
Исследовать, при каких $k$ разности этой функции при любом фиксированном $x$ неотрицательны
$$
\Delta^k g_n (x) > 0?
$$

Гипотеза: при любых $k$.

При $k=1$ гипотеза сводится к \eqref{43}, при $k=2$ к неравенству задачи 1 на стр. 20.

12. Верно ли, что при $n \in N_0$, $0<x<n+1$,
$$
R_{n+1, 1} (x) < R_{n, 1} (x)?
$$
К каким неравенствам для $R_n$ это приводит?

13. Исследовать отношения порядка между различными разностными аппроксимациями уравнения теплопроводности.

14. Перенести неравенства настоящей работы на остатки ряда Тейлора других простейших целых функций: $\ln (1+x)$, $\sin x$, $\cos x$, $J_{\nu} (x)$.

15. Найти $\max$, $\min$ и характер изменения при $x \in (0, \infty)$ функции
$$
f(x)= \frac{R_{n-2} R_n}{R^2_{n-1}} + \frac{R^2_{n}}{R_{n-1} R_{n+1}}.
$$
Из результатов настоящей работы следует, что
$$
1 + \frac{n}{n+1} = \frac{2 n+1}{n+1} \leq f(x) \leq 1 + \frac{n+2}{n+1} = 2 + \frac{1}{n+1} = \frac{2n+3}{n+1}.
$$
Отметим, что это очень узкие границы.

\break
\begin{center}
	Sitnik S.M. \\
	Inequalities for the exponential remainder of the Taylor series.\\
	Preprint, Vladivostok, 1993, \\
	AMS Subject Classification --- 26D15
\end{center}

\begin{center}
	Resume.
\end{center}

We study the remainder of the exponential function Taylor series defined by \eqref{1}. Recently Horst Alzer proved \cite{2} an inequality with the best possible constant \eqref{2}. The main aim of this preprint is consideration of \eqref{2} and its generalizations.

Th.1 states that \eqref{2} is equivalent to W.~Gautschi's inequality \eqref{6} from \cite{3}.

Th.2 is a simple modification, formally more general than \eqref{2}.

Th.3 contains reformations in terms of hypergeometric function, incomp\-lete gamma--function and fractional integrals.

Th.4 is a generalization of \eqref{2} again with the best constant. Proof is based on \eqref{19}--\eqref{20}, but \eqref{19} is simple and \eqref{20} is reduced to the Cebyshev inequality.

Another appropriate integral representation for the proof is $(18^*)$.

Introducing $(18^{**})$, we reduce \eqref{2} to \eqref{21}. It looks like interpolation\\ results.

Th.5 contains an interpolation inequality \eqref{22} with the best constant again. The proof is based on the Hadamard three--line theorem, \eqref{25}--\eqref{28} are consequences.

Reformulation of \eqref{2} as \eqref{29} looks like inequalities for absolutely and completely monotone functions but in different points.

Th.7 is a stronger result than \eqref{2}. The proof is reduced to an interesting inequality \eqref{32} which is trivial for integer $a=n$ but non--trivial for general $a$.

\eqref{36}--\eqref{40} are generalizations and consequences based on ideas of complete monotonicity. Kimberling's results froms \cite{22}, inequalities from \cite{19, 20, 21}.

The same results are true for more general term $R_{n, m} (x)$, occuring from the Obreshkov formula.

The interesting intriguing mistake is \eqref{42} with the wrong factorial. But if \eqref{42} were correct, then all the inequalities of this paper could be derived from the classical inequalities of Cebyshev, Jensen and Bunyakovsky by very simple procedures like \eqref{4}.

In the end there is a list of open problems and appendix with the reverse inequality \eqref{46} proved. It leads to the beautiful estimate \eqref{49} and connections with the Pade approximants.

\break

\end{document}